%% file: GGstrong.tex
\documentclass[12pt]{article}
\usepackage{amsmath}
\usepackage{amssymb}
\usepackage{vector}
\font\tencmmib=cmmib10 \skewchar\tencmmib '60
\newfam\cmmibfam
\textfont\cmmibfam=\tencmmib

\def\bbox{\quad\hbox{\vrule \vbox{\hrule \vskip2pt \hbox{\hskip2pt
\vbox{\hsize=1pt}\hskip2pt} \vskip2pt\hrule}\vrule}}
\def\lessim{\ \lower4pt\hbox{$
\buildrel{\displaystyle <}\over\sim$}\ }
\def\gessim{\ \lower4pt\hbox{$\buildrel{\displaystyle >}
\over\sim$}\ }

\def\la{{\Bigl\langle}}
\def\ra{{\Bigr\rangle}}

\def\qed{\hfill\break\rightline{$\bbox$}}
\newcommand{\EEE}{\mathbb{E}}

\newcommand{\Reals}{\mathbb{R}}

\newcommand{\vsi}{{\vec{\sigma}}}

\newtheorem{theorem}{Theorem}

\makeatletter
\@addtoreset{equation}{section}

\makeatother

\input invlat.tex

\begin{document}
\title{
The Ghirlanda-Guerra identities for mixed $p$-spin model.}

\author{ 
Dmitry Panchenko\thanks{
Department of Mathematics, Texas A\&M University, 
Mailstop 3386, College Station, TX, 77843,
email: panchenk@math.tamu.edu. Partially supported by NSF grant.}\\
{\it Texas A\&M University}
}
\date{}

\maketitle

\begin{abstract}
We show that under the conditions in \cite{T-P} which imply the validity of the Parisi formula
if the generic Sherrington-Kirkpatrick Hamiltonian contains a $p$-spin term then the 
Ghirlanda-Guerra identities for the $p$th power of the overlap hold in a strong sense without averaging. This implies strong version of the extended Ghirlanda-Guerra identities for mixed 
$p$-spin models than contain terms for all even $p\geq 2$ and $p=1.$
\end{abstract}
\vspace{0.5cm}



\section{Introduction and main result.}
A generic mixed $p$-spin Hamiltonian $H_N(\vsi)$ indexed by spin configurations 
$\vsi\in \{-1,+1\}^N$ is defined as a linear combination
\begin{equation}
H_N(\vsi) = \sum_{p\geq 1} \beta_p\, H_{p}(\vsi)
\label{HN}
\end{equation}
of $p$-spin Sherrington-Kirkpatrick Hamiltonians
\begin{equation}
H_{p}(\vsi)=\frac{1}{N^{(p-1)/2}}\sum_{1\leq i_1,\ldots,i_p \leq N}
g_{i_1,\ldots,i_p}\sigma_{i_1}\ldots \sigma_{i_p}
\label{pspin}
\end{equation}
where $(g_{i_1,\ldots,i_p})$ are i.i.d. standard Gaussian random variables, also independent
for all $p\geq 1.$ For simplicity of notation, we will keep the dependence of $H_p$ on $N$ implicit.
If a model involves an  external field parameter $h\in\Reals$ then the (random) Gibbs measure 
on $ \{-1,+1\}^N$ corresponding to the Hamiltonian $H_N$ is defined by
\begin{equation}
G_N(\vsi) = \frac{1}{Z_N} \exp\Bigl(H_N(\vsi)+h\sum_{i\leq N} \sigma_i\Bigr)
\end{equation}
where $Z_N$ is the normalizing factor called the partition function. As usual, we will denote by
$\la\cdot\ra$ the expectation with respect to the product Gibbs measure $G_N^{\otimes\infty}.$
One of the most important properties of the Gibbs measure $G_N$ was discovered by 
Ghirlanda and Guerra in \cite{GG} who showed that on average over some small perturbation
of the parameters $(\beta_p)$ in (\ref{HN}) the annealed product Gibbs measure satisfies a family
of distributional identities which are now called the Ghirlanda-Guerra identities. A more convenient
version of this result proved in \cite{Tal-New} can be formulated as follows. There exists a small perturbation $(\beta_{N,p})$ of the parameters in (\ref{HN}) such that all $\beta_{N,p}\to\beta_p$ 
and such that for all $p\geq 1$, $n\geq 2$ and any function
 $f=f(\vsi^1,\ldots,\vsi^n):(\{-1,+1\}^N)^n\to[-1,1]$, 
\begin{equation}
\lim_{N\to\infty}
\Bigl|
\EEE \la f  R_{1,n+1}^p \ra  - \frac{1}{n}\, \EEE \la f\ra \, \EEE\la R_{1,2}^p \ra -
\frac{1}{n}\sum_{l=2}^{n} \EEE \la f R_{1,l}^p\ra
\Bigr|
=0
\label{GGp}
\end{equation}
where $\la\cdot\ra$ is now the Gibbs average corresponding to the Hamiltonian (\ref{HN}) with
perturbed parameters $(\beta_{N,p})$ and $R_{l,l'}=N^{-1}\sum_{i\leq N}\sigma_i^l \sigma_i^{l'}$ 
is the overlap of configurations $\vsi^l$ and $\vsi^{l'}.$ Of course, the ultimate goal would be to show
that  (\ref{GGp}) holds without perturbing the parameters $(\beta_p)$ which would mean that the joint distribution of the overlaps $(R_{l,l'})_{l,l'\geq 1}$ under the annealed product Gibbs measure $\EEE G_N^{\otimes \infty}$ asymptotically satisfies the following distributional identities (up to symmetry considerations): for any $n\geq 2,$ conditionally on $(R_{l,l'})_{1\leq l<l' \leq n}$ the law of $R_{1,n+1}$ is given by the mixture $n^{-1} \mu + n^{-1} \sum_{l=2}^n \delta_{R_{1,l}}$ where $\mu$ is the law of $R_{1,2}.$
Toward this goal, a stronger version of (\ref{GGp}) for the original Hamiltonian (\ref{HN}) without
any perturbation of parameters was proved for $p=1$ in \cite{SCGG} under the additional assumption that $\beta_1\not = 0$ and a non-restrictive assumption on the limit of the free energy 
$F_N = N^{-1}\EEE\log Z_N.$ Here we will prove the same result for all $p$ under the assumptions 
and as a direct consequence of the seminal work of Talagrand in \cite{T-P} where the validity of the Parisi formula was proved. 
Namely,  from now on we will assume that the sum in (\ref{HN}) is taken over $p=1$ 
and even $p\geq 2.$ In this case, it was proved in \cite{T-P} that the limit of the free energy
\begin{equation}
\lim_{N\to\infty} F_N(\vec{\beta}) = P(\vec{\beta})
\label{Parisi}
\end{equation}
exists and is given by the Parisi formula $P(\vec{\beta})$ discovered in \cite{Parisi}. 
The exact definition of $P(\vec{\beta})$ will not be important to us and the only nontrivial property
that we will use is its differentiability in each coordinate $\beta_p$ which was
proved in \cite{PM} (see also \cite{PMP}).

\begin{theorem}\label{Th1}
For $p=1$ and for all even $p\geq 2,$
\begin{equation}
\lim_{N\to\infty} \frac{1}{N}\,\EEE 
\bigl\la
\bigl|
H_{p}(\vsi) - \EEE\la H_{p}(\vsi)\ra
\bigr|
\bigr\ra
=0.
\label{main}
\end{equation}
\end{theorem}
If $\beta_p\not = 0$ then (\ref{main}) implies (\ref{GGp}) by the usual integration by parts.
In particular, if $\beta_p \not =0$ for $p=1$ and all even $p\geq 2,$ the positivity principle 
of Talagrand proved in \cite{SG} implies the strong version of the extended Ghirlanda-Guerra
identities without any perturbation of the parameters. 


\textbf{Remark.} We will see that the proof  does not depend on the specific form of the Hamiltonian 
(\ref{HN}) and the result can be formulated in more generality.
Namely, given a sequence of random measures $\nu_N$ on some measurable space 
$(\Sigma,S)$ and a sequence of random processes $A_N$ indexed by $\vsi\in\Sigma$,
consider a sequence of Gibbs' measures $G_N$ defined by the change of density 
$$dG_N(\vsi)=Z_N^{-1}\exp (x A_N(\vsi)) d\nu_N(\vsi).$$
Let $\psi_N(x)=N^{-1}\log Z_N$ and $F_N(x)=\EEE \psi_N(x).$ Suppose that $\EEE|\psi_N(x)-F_N(x)|\to0$
and $F_N(x)\to P(x)$ in some neighborhood of $x_0$, and suppose that the limit $P(x)$ is differentiable at $x_0.$ Then
\begin{equation}
\lim_{N\to\infty} \frac{1}{N}\,\EEE 
\bigl\la
\bigl|
A_N(\vsi) - \EEE\la A_N(\vsi)\ra_{x_0}
\bigr|
\bigr\ra_{x_0}
=0,
\end{equation}
assuming some measurability and integrability conditions on $A_N$ and $\nu_N$
which will be clear from the proof and are usually trivially satisfied.
In Theorem \ref{Th1} we simply appeal to the results in \cite{PM} and \cite{T-P}.
\qed

\textbf{Proof of Theorem 1}. 
It has been observed for a long time that (see, for example, \cite{SCGG}) 
$$
\lim_{N\to\infty} \frac{1}{N}\,\EEE 
\bigl|
\la H_{p}(\vsi)\ra - \EEE\la H_{p}(\vsi)\ra
\bigr|
=0.
$$
This is where one uses the fact that $\EEE|\psi_N-F_N|\to0$, which is well known for 
$p$-spin models. It remains to prove that
\begin{equation}
\lim_{N\to\infty} \frac{1}{N}\,\EEE 
\bigl\la
\bigl|
H_{p}(\vsi) - \la H_{p}(\vsi)\ra
\bigr|
\bigr\ra
=0.
\end{equation}
This was proved in \cite{SCGG} for $p=1,$ but here we will show how this can be
obtained as a direct consequence of (\ref{Parisi}) for all $p$.
Let $\la\cdot\ra_x$ denote the Gibbs average corresponding to the Hamiltonian (\ref{HN}) 
where $\beta_p$ has been replaced by $x$. Consider $\beta_p^\prime>\beta_p$ and let 
$\delta=\beta_p^\prime - \beta_p.$ We start with the following obvious equation,
\begin{eqnarray}
\int_{\beta_p}^{\beta_p^\prime}\EEE\bigl\la
\bigl|
H_{p}(\vsi^1) - H_{p}(\vsi^2)
\bigr|
\bigr\ra_x
\,dx
&=&
\delta\, \EEE\bigl\la
\bigl|
H_{p}(\vsi^1) - H_{p}(\vsi^2)
\bigr|
\bigr\ra_{\beta_p}
\label{FTC}
\\
&&
+\,\int_{\beta_p}^{\beta_p^\prime}\! \! \!\int_{\beta_p}^{x}
\frac{\partial}{\partial t}
\EEE\bigl\la
\bigl|
H_{p}(\vsi^1) - H_{p}(\vsi^2)
\bigr|
\bigr\ra_t \,dt\, dx.
\nonumber
\end{eqnarray}
Since
\begin{eqnarray*}
\Bigl|
\frac{\partial}{\partial t}
\EEE\bigl\la
\bigl|
H_{p}(\vsi^1) - H_{p}(\vsi^2)
\bigr|
\bigr\ra_t 
\Bigr|
&=&
\Bigl|
\EEE\bigl\la
\bigl|
H_{p}(\vsi^1) - H_{p}(\vsi^2)
\bigr|
\bigl(H_{p}(\vsi^1)+H_{p}(\vsi^2)-2H_{p}(\vsi^3)\bigr)
\bigr\ra_t 
\Bigr|
\\
&\leq&
2\EEE\bigl\la
\bigl(
H_{p}(\vsi^1) - H_{p}(\vsi^2)
\bigr)^2
\bigr\ra_t 
\leq
8\EEE\bigl\la
\bigl(
H_{p}(\vsi) - \la H_{p}(\vsi) \ra_t
\bigr)^2
\bigr\ra_t 
\end{eqnarray*}
equation (\ref{FTC}) implies
\begin{eqnarray*}
\EEE\bigl\la
\bigl|
H_{p}(\vsi^1) - H_{p}(\vsi^2)
\bigr|
\bigr\ra_{\beta_p}
&\leq&
\frac{1}{\delta}
\int_{\beta_p}^{\beta_p^\prime}\EEE\bigl\la
\bigl|
H_{p}(\vsi^1) - H_{p}(\vsi^2)
\bigr|
\bigr\ra_x
\,dx
\\
&&
+\,
\frac{8}{\delta}
\int_{\beta_p}^{\beta_p^\prime}\! \! \!\int_{\beta_p}^{x}
\EEE\bigl\la
\bigl(
H_{p}(\vsi) - \la H_{p}(\vsi)\ra_t
\bigr)^2
\bigr\ra_t \,dt\, dx
\\
&\leq&
\frac{2}{\delta}
\int_{\beta_p}^{\beta_p^\prime}\EEE\bigl\la
\bigl|
H_{p}(\vsi) - \la H_{p}(\vsi)\ra_{x}
\bigr|
\bigr\ra_x
\,dx
\\
&&
+\,
8\int_{\beta_p}^{\beta_p^\prime}
\EEE\bigl\la
\bigl(
H_{p}(\vsi) - \la H_{p}(\vsi)\ra_t
\bigr)^2
\bigr\ra_t \,dt
\\
&\leq&
2\Bigl(\frac{1}{\delta}
\int_{\beta_p}^{\beta_p^\prime}\EEE\bigl\la
\bigl(
H_{p}(\vsi) - \la H_{p}(\vsi)\ra_{x}
\bigr)^2
\bigr\ra_x
\,dx\Bigr)^{1/2}
\\
&&
+\,
8\int_{\beta_p}^{\beta_p^\prime}
\EEE\bigl\la
\bigl(
H_{p}(\vsi) - \la H_{p}(\vsi)\ra_x
\bigr)^2
\bigr\ra_x \,dx.
\end{eqnarray*}
Therefore, if we denote
$$
\Delta_N = \frac{1}{N}
\int_{\beta_p}^{\beta_p^\prime}\EEE\bigl\la
\bigl(
H_{p}(\vsi) - \la H_{p}(\vsi)\ra_{x}
\bigr)^2
\bigr\ra_x
\,dx
$$
we showed that
\begin{equation}
\frac{1}{N}\,\EEE 
\bigl\la
\bigl|
H_{p}(\vsi) - \la H_{p}(\vsi)\ra
\bigr|
\bigr\ra
\leq
\frac{1}{N}\,\EEE 
\bigl\la
\bigl|
H_{p}(\vsi^1) - H_{p}(\vsi^2)
\bigr|
\bigr\ra
\leq
2\sqrt{\frac{\Delta_N}{N\delta}}
+8\Delta_N.
\label{Deltabound}
\end{equation}
If for a moment we think of $F_N=F_N(x)$ as a function of $x$ only then
$$
F_N'(x)=\frac{1}{N}\,\EEE\la H_p(\vsi)\ra_x \,\,\mbox{ and }\,\,
F_N''(x)
= \frac{1}{N}\,\EEE\bigl\la \bigl(
H_{p}(\vsi) - \la H_{p}(\vsi)\ra_{x}
\bigr)^2
\bigr\ra_x
$$
so that $\Delta_N = F_N'(\beta_p^\prime) - F_N^\prime(\beta_p).$
Since $F_N(x)$ is convex, for any $\gamma>0,$
$$
\Delta_N = F_N'(\beta_p^\prime) - F_N^\prime(\beta_p)
\leq
 \frac{F_N(\beta_p^\prime+\gamma) - F_N(\beta_p^\prime)}{\gamma}
-\frac{F_N(\beta_p) - F_N(\beta_p-\gamma)}{\gamma}
$$
and, therefore, equations (\ref{Deltabound}) and (\ref{Parisi}) now imply
$$
\limsup_{N\to\infty} \frac{1}{N}\EEE 
\bigl\la
\bigl|
H_{p}(\vsi) - \la H_{p}(\vsi)\ra
\bigr|
\bigr\ra
\leq
 8\Bigl(\frac{P(\beta_p^\prime+\gamma) - P(\beta_p^\prime)}{\gamma}
-\frac{P(\beta_p) - P(\beta_p-\gamma)}{\gamma}
\Bigr)
$$
where again we write $P=P(x)$ as a function of $x$ only.
Letting $\beta_p^\prime\to \beta_p$ first and then letting $\gamma\to 0$ and using that
$P(x)$ is differentiable proves the result.
\qed

\end{document}

%% file: invlat.tex
\font\tencmmib=cmmib10 \skewchar\tencmmib '60
\newfam\cmmibfam
\textfont\cmmibfam=\tencmmib

\def\bbox{\quad\hbox{\vrule \vbox{\hrule \vskip2pt \hbox{\hskip2pt
\vbox{\hsize=1pt}\hskip2pt} \vskip2pt\hrule}\vrule}}
\def\lessim{\ \lower4pt\hbox{$
\buildrel{\displaystyle <}\over\sim$}\ }
\def\gessim{\ \lower4pt\hbox{$\buildrel{\displaystyle >}
\over\sim$}\ }

\def\go0{\to 0}

\def\la{\langle}

\def\leftitem#1{\item{\hbox to\parindent{\enspace#1\hfill}}}

\def\qed{\hfill\break\rightline{$\bbox$}}

\def\ra{\rangle}

\def\sg{\sigma}

\def\sg2{\sigma^2}

\def\__{_{\infty}}

%% file: GGstrong.bbl
\begin{thebibliography}{99}

\bibitem{SCGG} S. Chatterjee, The Ghirlanda-Guerra identities without averaging,
preprint (2009), arXiv:0911.4520.

\bibitem{GG} S. Ghirlanda, F. Guerra, 
General properties of overlap probability distributions in disordered spin systems.
Towards Parisi ultrametricity,  J. Phys. A  31 (1998) no. 46 9149-9155.

\bibitem{PMP} D. Panchenko, On differentiability of the Parisi formula,
Elect. Comm. in Probab. 13 (2008) 241-247. 

\bibitem{Parisi} G. Parisi, A sequence of approximate solutions to the S-K model for spin glasses,
J. Phys. A 13 (1980) L-115.

\bibitem{PM} M. Talagrand, Parisi measures,  J. Funct. Anal. 231 (2006) no. 2 269-286.

\bibitem{T-P} M. Talagrand, Parisi formula, Ann. of Math. (2) 163 (2006) no. 1 221-263.

\bibitem{Tal-New} M. Talagrand, Construction of pure states in mean-field models for spin glasses, 
preprint (2008), to appear in Probab. Theory Related Fields.

\bibitem{SG} M. Talagrand, Mean Field Models for Spin Glasses. Manuscript.

\end{thebibliography}
